\documentclass{article}
\usepackage{graphicx} 
\usepackage{xcolor}
\usepackage{amsthm}
\usepackage{amsmath}
\usepackage{dsfont}
\usepackage{float}
\usepackage{natbib}
 \bibpunct[, ]{[}{]}{,}{n}{}{,}%
\usepackage{amssymb}

\newtheorem{theorem}{Theorem}[section]

\newtheorem{lemma}[theorem]{Lemma}
\newtheorem{remark}[theorem]{Remark}
\newtheorem{definition}[theorem]{Definition}

\newtheorem{def/prop}[theorem]{Definition/Proposition}

\newcommand{\Ltwo}{L^2([0,1])}
\newcommand{\model}{\mathbb{TIRG}}
\newcommand{\R}{\mathbf{R}}

\newcommand{\I}{\mathds{1}}

\newcommand{\procond}[2]{\mathbf{P}\mathopen{}\left(#1\middle|#2\right)}
\newcommand{\prob}[1]{\mathbf{P}\left( #1 \right)}
\newcommand{\NormOp}[1]{\left\vert\kern-0.25ex\left\vert\kern-0.25ex\left\vert #1 \right\vert\kern-0.25ex\right\vert\kern-0.25ex\right\vert}
\newcommand{\NormVect}[1]{\| #1 \|_{\R^{2^k},2 } }
\newcommand{\NormL}[1]{\| #1 \|_{L^2}}
\newcommand{\mat}[1]{\Lambda^{(#1)}}

\usepackage{enumitem}
\usepackage{authblk}

\SetLabelAlign{Center}{\hfil#1\hfil}
\SetLabelAlign{CenterWithParen}{(\makebox[1.0em]{#1})}

\title{Temporal Inhomogeneous Random Graphs}

\author[1,2]{Sofiya Burova\thanks{
Research of S.B. supported by the Spanish Agencia Estatal de Investigaci\'on under project PID2022-138268NB-100 and Grant PID2020-113082GB-I00 funded by MICIU/AEI/10.13039/501100011033.}}
\author[2,3]{Francisco Calvillo}
\date{\today}

\affil[1]{IMTECH and Departament de Matem\`atiques, Universitat Polit\`ecnica de Catalunya, Barcelona, Spain}
\affil[2]{Department of Economics and Business, Pompeu Fabra University, Barcelona, Spain}
\affil[3]{LPSM, Sorbonne Université, 4 Place Jussieu, 75005 Paris, France}

\begin{document}
\maketitle
\begin{abstract}
Motivated by the study of time-dependent propagation processes in inhomogeneous populations, we extend Burova, Lugosi and Perarnau's results on the temporal stochastic block model to graphon-based random graphs with a continuous and strictly positive kernel. We consider such a graph in the regime where connection probabilities are of the order of $\log n/n$, and equip its edges with a uniformly random ordering, interpreted as timestamps. We study the size of reachable sets via increasing paths (that is, paths whose edge timestamps are strictly increasing). We establish tight first-order asymptotics for the size of the reachable set from a fixed vertex, and also determine the first-order asymptotics of the lengths of the longest and shortest increasing paths with one or two fixed endpoints. 
\end{abstract}


\section{Introduction}
\subsection{Motivation and related works}

As most real-world networks evolve over time, propagation phenomena are more accurately described by models that incorporate temporal information. In many applications, such as epidemic spreading, information diffusion in social networks, or communication systems, the timing of interactions plays a crucial role: the existence of a connection between two vertices does not by itself guarantee that a signal, infection, or piece of information can propagate along it. Instead, propagation depends on the chronological order in which interactions occur. This motivates the study of graph models endowed with an explicit temporal structure.

A \emph{temporal graph} $G=(V,E,\pi)$ consists of a finite simple graph $(V,E)$ together with an ordering $\pi:E\to\{1,\ldots,|E|\}$ of its edges. The ordering $\pi$ is interpreted as assigning \emph{time stamps} to edges, representing the order in which interactions occur. An edge $e\in E$ is said to \emph{precede} an edge $f\in E$ if $$\pi(e)<\pi(f)~.$$
A path $(v_1,\ldots,v_k)$ from $u=v_1$ to $v=v_k$ is called an \emph{increasing path}, or \emph{temporal path}, if the time stamps of its consecutive edges are strictly increasing. In this case, $v$ is said to be \emph{reachable} from $u$. The \emph{reachable set} of a vertex $v\in V$ is the collection of all vertices that can be reached from $v$ through temporal paths. This induces a notion of \emph{temporal connectivity}, which is particularly relevant when studying propagation processes. For instance, if a vertex $v$ is initially infected by a virus, its reachable set represents the set of vertices that would become infected through the temporal evolution of the network.

\begin{figure}[H]
    \centering
    \includegraphics[scale=0.8]{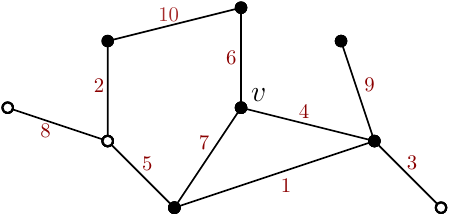}
    \caption{Example of a temporal graph. Time stamps are indicated next to the edges, and the vertices colored in black form the reachable set of vertex $v$.}
    \label{fig: reachable set}
\end{figure}

Motivated by the study of time-dependent propagation processes, \citet*{casteigts1}, and \citet*{casteigts2} introduced a natural random model of temporal graphs, called \emph{random simple temporal graphs}, obtained by equipping an Erd\H{o}s--R\'enyi random graph $G(n,p)$ with a uniformly random ordering of its edges. \citet{casteigts2} studied connectivity properties of the random simple temporal graph and showed that temporal reachability undergoes a phase transition at the characteristic scale
\[
p \asymp \frac{\log n}{n}.
\]

Further progress was made by \citet*{BrKaLu24}, who derived first-order asymptotics for the lengths of the longest and shortest temporal paths in this scale regime. Subsequent work by \citet*{AtDeLu25} refined the subcritical analysis of \citet{casteigts1}, showing that when $p$ lies below the threshold, the largest temporal clique (a set of vertices that are pairwise reachable from each other) has size of order $O(1)$ with high probability.

As a first step toward inhomogeneous temporal models, \citet*{sofiya} studied a temporal version of the stochastic block model, deriving first-order asymptotics for reachable sets and longest temporal paths. Their results show that the phase transition observed around $p=\log n / n$ for random simple temporal graphs persists in this more structured setting.

We continue this line of investigation by introducing a temporal version of
graphon-based random graphs, in which the finite collection of communities of
the stochastic block model is replaced by a continuum of vertex types.
Graphons were introduced by \citet{lovasz2006limits} and
\citet{BoChaLaSoVe12} as limiting objects of sequences of dense graphs. Their
sparse counterparts, known as \emph{inhomogeneous random graphs}, were studied
in depth by \citet*{BoJaRi07}, who showed that the emergence of a giant
component is governed by the norm of the integral operator associated with the
kernel. We show that the same operator governs temporal reachability, at the
scale $\log n/n$ instead of $1/n$. We start by recalling the definition of
graphon-based random graphs.

\subsection{Model and notation}
A \emph{graphon} is a symmetric Lebesgue-measurable function 
$f : [0,1]^2 \to [0,1]$.
Let $f$ be a graphon and let $n \ge 1$ be an integer. 
The random graph model on $n$ vertices associated with $f$, denoted by $G(n,f)$, is defined as follows. Let $(U_i)_{i\in [n]}$ be a sequence of i.i.d.\ random variables uniformly distributed on $[0,1]$. 
The vertex set is $V=[n]$, and, conditionally on $(U_i)_{i\in [n]}$, each edge $\{i,j\}$ is included independently with probability
\[
\procond{\{i,j\}\in E}{(U_k)_{k\in [n]}}
= f(U_i,U_j).
\]
\begin{definition}
Let $f:[0,1]^2\to[0,1]$ be a graphon and let $n\ge 1$. 
We define the \emph{Temporal Inhomogeneous Random Graph} model (TIRG for short), denoted by $\model(n,f)$, as a pair $(G,\pi)$, where
\begin{itemize}
    \item $G \sim G(n,f)$;
    \item $\pi : E(G) \to \{1,\ldots,|E(G)|\}$ is a uniformly random ordering of the edges of $G$, chosen independently of $G$.
\end{itemize}
\end{definition}

In other words, $\model(n,f)$ is the temporal random graph obtained by endowing the graphon-based random graph $G(n,f)$ with a uniform random ordering of its edges.

\begin{remark}
\label{unif_labels}
An equivalent construction of $\model(n,f)$ is obtained by
assigning to every edge of $G(n,f)$ an independent uniform
random variable on $[0,1]$, and ordering the edges according
to these values. Indeed, conditionally on the graph, all
orderings of its edges are equally likely.
\end{remark}

Observe that, when the graphon is constant, $f \equiv p \in [0,1]$, the graphon-based random graph $G(n,f)$ coincides with the classical Erd\H{o}s--R\'enyi random graph with parameters $(n,p)$. Consequently, the temporal model $\model(n,f)$ reduces to the random simple temporal graph studied in \cite{casteigts2,casteigts1,BrKaLu24,AtDeLu25}. This correspondence naturally suggests that temporal connectivity in the graphon-based setting should exhibit a phase transition around the scale $\log n / n$.

Throughout this work, we consider a sequence of graphons $(f_n)_{n\ge 1}$ of the form
\[
f_n = \rho_n \kappa,
\]
where
$\rho_n := \log n/n $ is the scaling parameter
and $\kappa:[0,1]^2 \to \mathbb{R}_{\ge 0}$ is a nonnegative, symmetric, bounded measurable function that does not depend on $n$. We implicitly assume that $n$ is sufficiently large so that $f_n$ takes values in $[0,1]$. Our interest lies in the associated sequence of the temporal inhomogeneous random graphs $\model(n,f_n)$.
More precisely, temporal connectivity is studied through the size of the reachable set from a fixed vertex $v$, denoted by
\[
B_n(v) := \{ u \in [n] : u \text{ is reachable from } v \}.
\]
A key quantity governing the propagation behavior in this model is the integral operator
$$K : \Ltwo \to \Ltwo$$ associated with the kernel $\kappa$, where $\Ltwo$
denotes the Hilbert space of real-valued square-integrable functions on $[0,1]$.
This operator plays a central role in the temporal connectivity behavior of $\model(n,f_n)$, in the same way as it governs the phase transition of the static inhomogeneous random graph with connection probabilities $\kappa(U_i,U_j)/n$; see \cite{BoJaRi07}.
For any $f \in \Ltwo$, the function $Kf$ is defined by
\begin{equation}
(Kf)(x) := \int_0^1 \kappa(x,y) f(y)\, dy
\ \text{ for any } x \in [0,1].
\end{equation}
We denote by $\NormOp{K}$ the operator norm of $K$, namely
\[
\NormOp{K}
= \sup_{f \neq 0}
\frac{\NormL{Kf}}{\NormL{f}},
\]
where $\NormL{\cdot}$ denotes the usual $L^2$-norm in $L^2([0,1])$.

Throughout the paper, all limits are taken as $n\to\infty$. We say that a
sequence of events $(E_n)_{n\ge1}$ holds \emph{with high probability} if
$\prob{E_n}\to1$. For a sequence of random variables $(X_n)_{n\ge1}$ and a
sequence of positive numbers $(a_n)_{n\ge1}$, we write $X_n\sim_p a_n$ if
$X_n/a_n\to1$ in probability. Finally, the \emph{length} of a path is its
number of edges, and we write $|\gamma|$ for the length of a path $\gamma$.

\subsection{Main results}
The main contribution of this paper is summarized in the following theorem. 
It reveals a phase transition phenomenon for temporal connectivity around the critical value $\NormOp{K}=1$ and provides first-order asymptotics for the size of the reachable set of a fixed vertex.

\begin{theorem}\label{thm: phase transition}
    Let $\kappa:[0,1]^2\to\R_{\ge0}$ be a symmetric, bounded continuous
function such that
\[
\min_{(x,y)\in[0,1]^2}\kappa(x,y)>0~,
\] and let $\rho_n = \frac{\log n }{n}$ for any $n\ge 1~.$ Then, for any fixed vertex $v$ in $\model(n,\rho_n \kappa)~$,
    \begin{enumerate}
    \item[(i)]\label{i}if $\NormOp{K}>1$, then $$|B_n(v)| = (1-o(1))n$$ with high probability;
    \item[(ii)] if $\NormOp{K}<1$, then for all $\epsilon >0$, 
    $$n^{\NormOp{K}(1-\epsilon)}\le |B_n(v)|\le n^{\NormOp{K}(1+\epsilon)}$$
    with high probability.
    \end{enumerate}
\end{theorem}

\medskip

In addition to the size of the reachable set, we study the lengths of
increasing paths. Let $\gamma_{\max}$, $\gamma_{\max}(v)$, and $\gamma_{\max}(u,v)$ denote
the longest increasing path, the longest increasing path starting at
$v$, and the longest increasing path from $u$ to $v$, respectively.
Similarly, let $\gamma_{\min}(u,v)$ denote the shortest increasing path
from $u$ to $v$, whenever such a path exists.

\begin{theorem}\label{thm: graphon longest paths}
Let $\kappa:[0,1]^2\to\R_{\ge0}$ be symmetric, continuous and bounded
away from zero. Write $\theta=\NormOp{K}$ and fix two distinct vertices
$u,v\in[n]$. Then,
\begin{enumerate}
\item[(i)] for every $\epsilon>0$, with high probability,
\[
(e\theta-\epsilon)\log n
\le |\gamma_{\max}|
\le(\alpha(\theta)+\epsilon)\log n~,
\]
where
\[
\alpha(\theta):=\inf\left\{x>0:
 x\log\left(\frac{x}{e\theta}\right)=1\right\}~;
\]
\item[(ii)] $|\gamma_{\max}(v)|\sim_p e\theta\log n$;
\item[(iii)] if $\theta>1$, then
$|\gamma_{\max}(u,v)|\sim_p\beta(\theta)\log n$ and
$|\gamma_{\min}(u,v)|\sim_p\gamma(\theta)\log n$, where
\[
\beta(\theta):=\sup\left\{x>0:
 x\log\left(\frac{x}{e\theta}\right)=-1\right\}~,
\]
\[
\gamma(\theta):=\inf\left\{x>0:
 x\log\left(\frac{x}{e\theta}\right)=-1\right\}~.
\]
If $\theta<1$, then, with high probability, there is no increasing
path from $u$ to $v$.
\end{enumerate}
\end{theorem}

\section{Stochastic Block Model Approximations}\label{sec: sbm approximations}

To obtain lower and upper bounds for the reachable set, we approximate the graphons
$(f_n)_{n\ge1}$ by step functions. Such approximations partition the unit interval
into finitely many blocks and replace the kernel by a piecewise constant function.
The associated random graphs are stochastic block models, where vertices are grouped
into finitely many types and connection probabilities depend only on the types. This connection allows us to relate the TIRG model to temporal stochastic
block models, for which results on reachable sets are already available; see Theorem~1.1 in \cite{sofiya}. The approximation will therefore be used to transfer these results
and obtain lower and upper bounds in both the subcritical and supercritical regimes.

\medskip

Let $k\ge 1$ be an integer. For each $i\in [2^k]$, denote by
\[
J_i^k := \Big[\frac{i-1}{2^k},\frac{i}{2^k}\Big)
\]
the dyadic partition of $[0,1]$. Define the function
$\kappa_k : [0,1]^2 \to \R_{\ge 0}$ by
\[
\forall i,j\in [2^k],\ \forall x\in J_i^k,\ \forall y\in J_j^k,
\qquad
\kappa_k(x,y)
=
\min_{\overline{J_i^k}\times\overline{J_j^k}}\kappa .
\]
Similarly, define the upper step approximation $\overline\kappa_k$ by
\[
\forall i,j\in[2^k],\ \forall x\in J_i^k,\ \forall y\in J_j^k,
\qquad
\overline\kappa_k(x,y)
=
\max_{\overline{J_i^k}\times\overline{J_j^k}}\kappa .
\]
The minimum and maximum are attained on the closures of the dyadic
rectangles, since $\kappa$ is continuous. In particular,
\[
\kappa_k\le\kappa\le\overline\kappa_k~.
\]
Let $K_k,\overline K_k:L^2([0,1])\to L^2([0,1])$ be the integral
operators associated with $\kappa_k$ and $\overline\kappa_k$, respectively,
defined for any $f\in L^2([0,1])$ by
\[
(K_k f)(x) = \int_0^1 f(y)\kappa_k(x,y)\,dy,
\qquad
(\overline K_k f)(x) = \int_0^1 f(y)\overline\kappa_k(x,y)\,dy .
\]
All the results below apply to both approximations, with minima replaced
by maxima for the upper approximation. The direction of the graph
inclusion is reversed in the latter case.
Since $\kappa_k$ is a step function, the random graph
$G(n,\rho_n\kappa_k)$ follows a stochastic block model with $2^k$
communities, conditionally on the sequence $(U_i)_{i\ge 1}$. More precisely, each vertex $j\in[n]$ is assigned to the
community indexed by the interval $J_i^k$ containing the latent variable
$U_j$, and edges between communities $i$ and $j$ are present with
probability
\[
p_{ij}=\frac{\log n}{n}\min_{\overline{J_i^k}\times\overline{J_j^k}}\kappa .
\]
For the upper step model, the edge probabilities are obtained by
replacing the minimum by the maximum.
The community partition is given by the sets 
\[
V_i^{k } = \{j\in [n]: U_j\in J_i^k\}~.
\]
Since the variables $(U_j)_{j\ge1}$ are i.i.d.\ uniform, each vertex is assigned to community $i$ independently with probability $2^{-k}$. Hence $\model(n,\rho_n\kappa_k)$ is exactly the temporal stochastic block model with \emph{random blocks} described in \cite[Section~5.1]{sofiya}, with $2^k$ communities, proportions $(2^{-k},\ldots,2^{-k})$ and connection probabilities $p_{ij}$. The model for which the results of \cite{sofiya} are stated has deterministic blocks instead, of sizes $\lfloor n/2^k\rfloor$ or $\lceil n/2^k\rceil$. As observed in \cite[Section~5.1]{sofiya}, the two models are asymptotically equivalent and the results of that paper hold for the random-block version as well. We rely on this fact throughout: Theorems~1.1 and~1.6 of \cite{sofiya} will be applied to $\model(n,\rho_n\kappa_k)$ and $\model(n,\rho_n\overline\kappa_k)$. The main ingredient behind this equivalence is the law of large numbers, which we record in the next lemma.
\begin{lemma}\label{lemma: LLN sbm}
Let $k\ge 1$ be fixed. Then, with high probability,
\[
(|V_1^{k}|,\ldots,|V_{2^k}^{k}|)
=
\left(\frac{n}{2^k},\ldots,\frac{n}{2^k}\right)
+ o(n).
\]
\end{lemma}

\begin{proof}
For each fixed $i\in[2^k]$,
\[
|V_i^{k}|
=
\sum_{j=1}^n \mathbf{1}_{\{U_j\in J_i^k\}}.
\]
Since $(U_j)_{j\ge1}$ are i.i.d.\ uniform on $[0,1]$, we have
\[
\prob{U_j\in J_i^k}=\frac{1}{2^k}.
\]
Hence $|V_i^{k}|$ is a binomial random variable of parameters $n$ and $1/2^k$, so by the law of large numbers,
\[
\frac{|V_i^{k}|}{n}\longrightarrow \frac{1}{2^k}
\quad \text{almost surely}.
\]
Since $2^k$ is fixed, a union bound over $i\in[2^k]$ yields the result.
\end{proof}
For $n,k\ge1$ and a vertex $v\in[n]$, we denote by $B_n^k(v)$
and $\overline B_n^k(v)$ the reachable sets of $v$ in the temporal
graphs $\model(n,\rho_n\kappa_k)$ and
$\model(n,\rho_n\overline\kappa_k)$, respectively.
\begin{lemma}\label{lemma: sbm subgraph}
For every $n,k\ge1$, there exists a coupling of the three temporal
graphs under which
\[
G(n,\rho_n\kappa_k)
\subseteq
G(n,f_n)
\subseteq
G(n,\rho_n\overline\kappa_k).
\]
Consequently, for every vertex $v$ and every $\lambda\in\R$,
\begin{equation}\label{eq: sbm vs graphon}
\prob{|\overline B_n^k(v)|\le\lambda}
\le
\prob{|B_n(v)|\le \lambda}
\le
\prob{|B_n^k(v)|\le \lambda}.
\end{equation}
\end{lemma}

\begin{proof}
Since $\kappa_k(x,y)\le\kappa(x,y)\le\overline\kappa_k(x,y)$
for all $(x,y)\in[0,1]^2$, we have
\[
\rho_n\kappa_k(x,y)\le f_n(x,y)
\le\rho_n\overline\kappa_k(x,y)~.
\]
We construct the three graphs using the same latent variables
$(U_i)_{i\ge1}$ and two further families of i.i.d.\ uniform random
variables $(\xi_{ij})_{1\le i<j\le n}$ and $(\tau_{ij})_{1\le i<j\le n}$,
the three families being independent. An edge $\{i,j\}$ is present
in each graph when $\xi_{ij}$ does not exceed the corresponding edge
probability, and in that case it receives the time stamp $\tau_{ij}$ in
every graph containing it. By Remark~\ref{unif_labels}, each of the three
temporal graphs has the correct distribution, and the stated inclusions
hold almost surely.
Temporal paths are preserved under edge inclusion, therefore
\[
B_n^k(v)\subseteq B_n(v)\subseteq\overline B_n^k(v)
\]
under this coupling, which implies \eqref{eq: sbm vs graphon}.
\end{proof}

As proven in \cite{sofiya}, the temporal connectivity (and its phase transition) in the stochastic block model is driven by the spectral radius of a matrix that can be seen as the finite-dimensional analogue of the operator $K$. To make the connection explicit, we introduce the matrix $\mat{k}\in \R^{2^k\times2^k}$ whose coordinates are $$\mat{k}_{ij}=2^{-k} \min_{\overline{J_i^k}\times\overline{J_j^k}} \kappa$$ and let $\theta_k$ be its spectral radius, i.e.
\[
\theta_k := \max \left\{|\lambda|: \lambda \in \sigma (\Lambda^{(k)})\right\}~,
\]
where $\sigma (\Lambda^{(k)})$ is the set of eigenvalues of $\Lambda^{(k)}~.$
Since $\mat{k}$ is symmetric, $\theta_k$ is equal to the spectral norm of $\mat{k}$, i.e.
\begin{equation}
    \theta_k = \sup_{u\in \R^{2^k}, u\neq 0} \frac{\NormVect{\mat{k}u}}{\NormVect{u}}~, \label{eqSpectralNorm}
\end{equation}
where $\NormVect{\cdot}$ denotes the classic Euclidean norm in $\R^{2^k}~.$
Similarly, let $\overline\Lambda^{(k)}$ be the matrix with coordinates
\[
\overline\Lambda^{(k)}_{ij}
=2^{-k}\max_{\overline{J_i^k}\times\overline{J_j^k}}\kappa
\]
and let $\overline\theta_k$ be its spectral radius.

\begin{lemma}\label{lemma: eigenvalue equals norm}
    $\theta_k = \NormOp{K_k}$ and
    $\overline\theta_k = \NormOp{\overline K_k}$.
\end{lemma}

\begin{proof}
Fix $k$ and observe that for any $f\in L^2([0,1])$ and any $x\in J_i^k$ for some fixed $i\in [2^k]$,
\begin{align*}
    K_k f(x) &= \int_0^1 \kappa_k(x,y)f(y)dy\\
    &=\sum_{j=1}^{2^k} \int_{\frac{j-1}{2^k}}^\frac{j}{2^k} \kappa_k(x,y)f(y)dy \\
    &= 2^{k}\sum_{j=1}^{2^k} \int_{\frac{j-1}{2^k}}^\frac{j}{2^k} \mat{k}_{ij}f(y)dy~.
\end{align*}
Denote by $v_f$ the vector in $\R^{2^k}$ whose coordinates are $(v_f)_j = \int_{J^k_j}f(y)dy$. Then $K_k f (x) = 2^k\sum_{j=1}^{2^k}  \mat{k}_{ij}(v_f)_j$ whenever $x\in J_i^k$.
Thus,
\begin{align}
    \NormL{K_k f}^2 &= \sum_{i=1}^{2^k} \int_{J_i^k} \left( 2^{k}\sum_{j=1}^{2^k} \mat{k}_{ij}(v_f)_j \right)^2 dx\nonumber \\
    &= \frac{1}{2^k}\sum_{i=1}^{2^k} \left( 2^k\sum_{j=1}^{2^k} \mat{k}_{ij}(v_f)_j \right)^2 \nonumber \\
    &= 2^k\NormVect{\mat{k}v_f}^2~. \label{eqNormes}
\end{align}
From \eqref{eqSpectralNorm} and \eqref{eqNormes} it follows that for any $f\in L^2([0,1])$,
\begin{equation}
    \NormL{K_k f}^2 \le 2^k \theta_k^2 \NormVect{v_f}^2~. \label{eqOperatorToMatrix}
\end{equation}
Moreover, by the Cauchy--Schwarz inequality, we have $(v_f)_i^2 \le 2^{-k}\int_{J_i^k} f(x)^2dx$, leading to $$\NormVect{v_f}^2\le 2^{-k} \NormL{f}^2~.$$
Together with \eqref{eqOperatorToMatrix}, we have that for any $f\in L^2([0,1])$, $$\NormL{K_k f}^2\le\theta_k^2 \NormL{f}^2~,$$ which implies that $$\NormOp{K_k}\le \theta_k~.$$
Now, let $u\in \R^{2^k}$ be such that $\NormVect{u}=1$ and $\mat{k}u= \theta_k u$. Define $f_{\theta_k} \in L^2([0,1])$ by $f_{\theta_k}(x)=2^k \sum_{i=1}^{2^k}\I_{x\in J_i^k} u_i$, so that we have $ v_{f_{\theta_k}}=u$. Its norm is $$\NormL{f_{\theta_k}}^2=4^k\sum_{i=1}^{2^k}2^{-k} u_i^2 =2^k~.$$
Using equality \eqref{eqNormes} it follows that 
\begin{align*}
    \theta_k ^2 &= \NormVect{\mat{k}u}^2 \\
    &=2^{-k}\NormL{K_k f_{\theta_k}}^2 \\
    &\le 2^{-k}\NormOp{K_k}^2 \NormL{f_{\theta_k}}^2 = \NormOp{K_k}^2~,
\end{align*}
proving the first equality. The same computation, with maxima in place
of minima, gives $\overline\theta_k=\NormOp{\overline K_k}$~.
\end{proof}

\begin{lemma}\label{lemma: operator convergence}
Assume that $\kappa$ is a nonnegative, symmetric, bounded continuous function. Then
    \[
    \NormOp{K - K_k}\to 0
    \quad\text{and}\quad
    \NormOp{\overline K_k-K}\to 0
    \quad \text{ as } \quad k \to \infty~.
    \]
\end{lemma}
\begin{proof}
Since $\kappa$ is continuous on the compact set $[0,1]^2$, it is uniformly
continuous. Let
\[
\omega_\kappa(\delta)
:=
\sup\left\{
    |\kappa(x,y)-\kappa(x',y')| :
    \left\|(x,y)-(x',y')\right\|_2\le \delta
\right\}
\]
denote its modulus of continuity. Then
\[
\omega_\kappa(\delta)\to 0
\qquad \text{as } \quad \delta\to 0~.
\]
For any $1\le i,j\le2^k$ the diameter of $J_i^k\times J_j^k$ is
$\sqrt{2}\,2^{-k}$. Hence
\[
0
\le
\kappa(x,y)-\kappa_k(x,y)
\le
\omega_\kappa\left(\sqrt{2}\,2^{-k}\right),
\qquad
0\le\overline\kappa_k(x,y)-\kappa(x,y)
\le\omega_\kappa\left(\sqrt{2}\,2^{-k}\right)
\]
for any $x,y \in [0,1]~.$
Taking the supremum over $(x,y)\in[0,1]^2$, we obtain
\[
\sup_{(x,y)\in[0,1]^2}
|\kappa(x,y)-\kappa_k(x,y)|
\le
\omega_\kappa\left(\sqrt{2}\,2^{-k}\right)
\longrightarrow 0~.
\]
Finally, for every $f\in L^2([0,1])$,
\begin{align*}
\NormL{(K-K_k)f}^2
&=
\int_0^1
\left|
\int_0^1
\bigl(\kappa(x,y)-\kappa_k(x,y)\bigr)f(y)\,dy
\right|^2
dx \\
&\le
\left(
\sup_{(x,y)\in[0,1]^2}
|\kappa(x,y)-\kappa_k(x,y)|
\right)^2
\NormL{f}^2~,
\end{align*}
where we used the Cauchy--Schwarz inequality and the fact that
$[0,1]$ has Lebesgue measure one. Therefore,
\[
\NormOp{K-K_k}
\le
\sup_{(x,y)\in[0,1]^2}
|\kappa(x,y)-\kappa_k(x,y)|
\longrightarrow 0~,
\]
and the same argument gives $\NormOp{\overline K_k-K}\to0$,
which proves the result.
\end{proof}
\section{Reachable sets}
\subsection{Supercritical regime}

\begin{proof}[Proof of Theorem~\ref{thm: phase transition}(i)]
Assume that $\NormOp{K}>1$ and fix $\epsilon>0$.

\noindent By Lemma~\ref{lemma: sbm subgraph}, for every $k$,
\begin{equation}\label{eq: supercritical inequality}
\prob{\frac{|B_n(v)|}{n} \le 1-\epsilon}
\le
\prob{\frac{|B_n^k(v)|}{n}\le 1-\epsilon}.
\end{equation}
Since $\NormOp{K}>1$, there exists $\eta>0$ such that
\[
(1-\eta)\NormOp{K}>1.
\]
By Lemmas~\ref{lemma: eigenvalue equals norm} and
\ref{lemma: operator convergence}, we may choose $k$ large enough so that
\[
\theta_k \ge (1-\eta)\NormOp{K}>1.
\]
Hence the finite-dimensional matrix governing the stochastic block
model associated with $\kappa_k$ has spectral radius strictly larger
than $1$, and the model lies in the supercritical regime studied
in \cite{sofiya}. Since $\kappa$ is bounded away from zero, every entry of
$\mat{k}$ is strictly positive. In particular, $\mat{k}$ is
irreducible, so the irreducibility assumption of Theorem~1.1
in \cite{sofiya} is satisfied.
As explained in Section~\ref{sec: sbm approximations}, Theorem~1.1(i) of
\cite{sofiya} applies to $\model(n,\rho_n\kappa_k)$ and gives
$|B_n^k(v)|\sim_p n$. In particular,
\[
\prob{\frac{|B_n^k(v)|}{n}\le 1-\epsilon}
=
o(1)
\quad \text{as } n\to\infty.
\]
Combining this with
\eqref{eq: supercritical inequality}, we conclude that
\[
\prob{\frac{|B_n(v)|}{n} \le 1-\epsilon}
=
o(1)
\quad \text{as } n\to\infty,
\]
which completes the proof.
\end{proof}

\subsection{Subcritical regime}
\begin{proof}[Proof of Theorem~\ref{thm: phase transition}(ii)]
Assume that $\NormOp{K}<1$ and fix $\epsilon>0$.
We first establish the lower bound. By Lemmas~\ref{lemma: eigenvalue equals norm} and~\ref{lemma: operator convergence}, there exists $k\ge1$ such that
\[
\NormOp{K} \le (1+\epsilon)\theta_k.
\]
It follows that
\begin{align}
\prob{|B_n(v)|\le n^{(1-\epsilon)\NormOp{K}}}
&\le \prob{|B_n(v)|\le n^{(1-\epsilon^2)\theta_k}} \nonumber\\
&\le \prob{|B_n^k(v)|\le n^{(1-\epsilon^2)\theta_k}}, 
\label{eq: subcritical lowebound}
\end{align}
where the last inequality follows from Lemma~\ref{lemma: sbm subgraph}.

Since $\theta_k\to\NormOp{K}<1$, taking $k$ larger if necessary we may assume that $\theta_k<1$, so the stochastic block model associated with $\kappa_k$ lies in the subcritical regime studied in \cite{sofiya}. Moreover, $\mat{k}$ is irreducible. As explained in Section~\ref{sec: sbm approximations}, Theorem~1.1(ii) of \cite{sofiya} applies to $\model(n,\rho_n\kappa_k)$. Applying it with $\epsilon^2/2$ in place of $\epsilon$, we get
\[
|B_n^k(v)| \ge n^{(1-\epsilon^2/2)\theta_k} > n^{(1-\epsilon^2)\theta_k}
\]
with high probability. Combining this with \eqref{eq: subcritical lowebound} yields
\[
\prob{|B_n(v)| \le n^{(1-\epsilon)\NormOp{K}}}
= o(1)
\quad \text{as } n\to\infty.
\]
The upper bound follows by using $\overline\kappa_k$ instead of
$\kappa_k$. Indeed, by Lemmas~\ref{lemma: eigenvalue equals norm}
and~\ref{lemma: operator convergence},
$\overline\theta_k\to\NormOp{K}$, and by
Lemma~\ref{lemma: sbm subgraph},
$B_n(v)\subseteq\overline B_n^k(v)$ under the coupling. Choose $k$ sufficiently large that $\overline\theta_k<1$ and
$(1+\eta)\overline\theta_k<(1+\epsilon)\NormOp{K}$ for some
$\eta>0$. The upper bound from Theorem~1.1 in \cite{sofiya}, applied
to the upper step model, then gives
$|B_n(v)|\le n^{(1+\epsilon)\NormOp{K}}$ with high probability.

This completes the proof.
\end{proof}

\section{Longest paths}

\begin{proof}[Proof of Theorem~\ref{thm: graphon longest paths}]
Fix $\epsilon>0$. Let $\kappa_k$ and $\overline\kappa_k$ be the lower and upper
step approximations introduced in Section~\ref{sec: sbm approximations}.
Write $\theta_k=\NormOp{K_k}$ and
$\overline\theta_k=\NormOp{\overline K_k}$.
By Lemmas~\ref{lemma: eigenvalue equals norm}
and~\ref{lemma: operator convergence},
\[
\theta_k,\overline\theta_k\longrightarrow\NormOp{K}~.
\]
We construct the three temporal graphs using the same latent
variables, edge-presence variables and timestamps (as defined in Remark \ref{unif_labels}). Then,
\[
G(n,\rho_n\kappa_k)\subseteq G(n,\rho_n\kappa)
\subseteq G(n,\rho_n\overline\kappa_k)~,
\]
and the same inclusion holds for their increasing paths.
For fixed $k$, the lower and upper step models are temporal stochastic
block models with random blocks, with $2^k$ communities of proportions
$2^{-k}$; see Section~\ref{sec: sbm approximations}. Since $\kappa$ is
bounded away from zero, their parameter matrices are irreducible. By
\cite[Section~5.1]{sofiya}, the path-length estimates of Theorem~1.6 in
\cite{sofiya} therefore apply to both models.

Let $\gamma_{\max}^k(v)$ and $\overline\gamma_{\max}^k(v)$ denote
the longest increasing paths from $v$ in the lower and upper step
models. Under the above coupling,
\[
|\gamma_{\max}^k(v)|\le|\gamma_{\max}(v)|
\le|\overline\gamma_{\max}^k(v)|~.
\]
By Theorem~1.6(ii) of \cite{sofiya}, the first and last quantities
are asymptotic to $e\theta_k\log n$ and
$e\overline\theta_k\log n$, respectively. Letting $n\to\infty$
and then $k\to\infty$ proves (ii). The lower bound in (i) follows
from (ii), while the upper bound follows by applying
Theorem~1.6(i) of \cite{sofiya} to the upper step model and
letting $k\to\infty$, using the continuity of $\alpha$.

Suppose next that $\NormOp{K}>1$. Choose $k$ large enough that
$\theta_k>1$. The lengths of the longest increasing paths between u and v are ordered in the same way. For the shortest paths, the inequalities are
reversed, namely,
\[
|\overline\gamma_{\min}^k(u,v)|
\le|\gamma_{\min}(u,v)|
\le|\gamma_{\min}^k(u,v)|~,
\]
with high probability, since the shortest path exists in the lower
step model. Theorem~1.6(iii) of \cite{sofiya} gives the corresponding
constants $\beta(\theta_k),\beta(\overline\theta_k)$ for the
longest paths and $\gamma(\theta_k),\gamma(\overline\theta_k)$
for the shortest paths. By continuity of $\beta$ and $\gamma$ for
$\theta>1$, letting $k\to\infty$ proves (iii) in this regime.

Finally, suppose $\NormOp{K}<1$. Choose $k$ large enough that
$\overline\theta_k<1$. Theorem~1.6(iii) of \cite{sofiya} shows
that, with high probability, there is no increasing path from $u$
to $v$ in the upper step model. The same holds in
$\model(n,\rho_n\kappa)$, completing the proof.
\end{proof}

\bibliographystyle{plainnat}
\bibliography{biblio.bib}

@article{casteigts1,
  title={Giant components in random temporal graphs},
  author={Becker, Ruben and Casteigts, Arnaud and Crescenzi, Pierluigi and Kodric, Bojana and Renken, Malte and Raskin, Michael and Zamaraev, Viktor},
  journal={arXiv preprint arXiv:2205.14888},
  year={2022}
}

@article{casteigts2,
  title={Sharp thresholds in random simple temporal graphs},
  author={Casteigts, Arnaud and Raskin, Michael and Renken, Malte and Zamaraev, Viktor},
  journal={SIAM Journal on Computing},
  volume={53},
  number={2},
  pages={346--388},
  year={2024},
  publisher={SIAM}
}

@article{BrKaLu24,
  title={Increasing paths in random temporal graphs},
  author={Broutin, Nicolas and Kam{\v{c}}ev, Nina and Lugosi, G{\'a}bor},
  journal={The Annals of Applied Probability},
  volume={34},
  number={6},
  pages={5498--5521},
  year={2024},
  publisher={Institute of Mathematical Statistics}
}

@article{AtDeLu25,
  title={On the size of temporal cliques in subcritical random temporal graphs},
  author={Atamanchuk, Caelan and Devroye, Luc and Lugosi, G{\'a}bor},
  journal={Combinatorics, Probability and Computing},
  volume={34},
  number={5},
  pages={671--679},
  year={2025},
  publisher={Cambridge University Press}
}

@article{lovasz2006limits,
  title={Limits of dense graph sequences},
  author={Lov{\'a}sz, L{\'a}szl{\'o} and Szegedy, Bal{\'a}zs},
  journal={Journal of Combinatorial Theory, Series B},
  volume={96},
  number={6},
  pages={933--957},
  year={2006},
  publisher={Elsevier}
}

@article{BoChaLaSoVe12,
author = {Borgs, Christian and Chayes, Jennifer and Lovász, László and Sós, Vera and Vesztergombi, Katalin},
year = {2012},
month = {07},
pages = {},
title = {Convergent sequences of dense graphs. II. Multiway cuts and statistical physics},
volume = {1},
journal = {Annals of Mathematics. Second Series},
doi = {10.4007/annals.2012.176.1.2}
}

@article{sofiya,
      title={The temporal stochastic block model}, 
      author={Burova, S. and Lugosi, G. and Perarnau, G.},
      year={2026},
  journal = {arXiv preprint arXiv:2606.09488}
}

@article{BoJaRi07,
  author  = {Bollob{\'a}s, B{\'e}la and Janson, Svante and Riordan, Oliver},
  title   = {The phase transition in inhomogeneous random graphs},
  journal = {Random Structures \& Algorithms},
  volume  = {31},
  number  = {1},
  pages   = {3--122},
  year    = {2007}
}
\end{document}